\documentclass[a4paper,11pt]{article}   
\usepackage{lineno,hyperref}
\usepackage[utf8]{inputenc}
\usepackage[english]{babel}
\usepackage{ifthen}
\usepackage{amssymb,amsmath}
\usepackage[usenames]{color}
\date{}

 \newif\ifNoRemark
    \def\addtheorem#1#2#3#4{ 
    \ifthenelse{\expandafter\isundefined\csname the#2\endcsname}{\newcounter{#2}}{}
    \newenvironment{#1}[1][\global\NoRemarktrue]
     {\par\addvspace{2mm}\noindent
       \refstepcounter{#2}{\bf #3~\csname the#2\endcsname
      \vphantom{##1}\ifNoRemark.\ \else\ (##1).\fi}\begingroup #4}%
     {\endgroup\par\addvspace{1mm}\global\NoRemarkfalse}
    \expandafter\newcommand\csname b#1\endcsname{\begin{#1}}
    \expandafter\newcommand\csname e#1\endcsname{\end{#1}}
    }

\addtheorem{theorem}{thrm}{Theorem}{\sl}
\addtheorem{lemma}{lmm}{Lemma}{\sl}
\addtheorem{proposition}{prpstn}{Proposition}{\sl}
\addtheorem{corollary}{crllr}{Corollary}{\sl}
\addtheorem{problem}{problem}{Problem}{\sl}
\addtheorem{remark}{}{Remark}{}

\begin{document}

\title{DP-colorings of uniform hypergraphs and splittings of
Boolean hypercube into faces
\thanks{The work was funded by the Russian Science
Foundation (grant No 18-11-00136).}}

\author{ Vladimir N. Potapov\\ {\em \small Sobolev Institute of Mathematics, Novosibirsk,
Russia; email: vpotapov@math.nsc.ru}}

 \maketitle

\begin{abstract}

We develop a connection between   DP-colorings  of    $k$-uniform
hypergraphs of order $n$ and coverings of $n$-dimensional Boolean
hypercube by  pairs of antipodal $(n-k)$-dimensional faces.
Bernshteyn  and Kostochka established that the lower bound on edges
in a non-2-DP-colorable $k$-uniform
 hypergraph is equal to $2^{k-1}$ for odd $k$ and $2^{k-1}+1$ for even $k$.
 They proved that these bounds are tight for $k=3,4$. In this paper,  we prove
 that the bound is achieved for all odd $k\geq 3$.
 \end{abstract}


MSC  05C15, 05C65, 05C35,  51E05

\section{Introduction}

Let $Q_2^n$ be an $n$-dimensional  Boolean hypercube.  We consider
coverings and splittings of $Q_2^n$ into faces.  A $k$-coverings of
$Q_2^n$ is a set of $(n-k)$-dimensional axis-aligned planes or
$(n-k)$-faces such that a union of the faces is equal to $Q_2^n$.
Two $m$-faces are called parallel if they have the same directions
and a pair of parallel faces is called antipodal if for each vertex
from one face there exists an antipodal vertex in another face.  It
is clear that each $k$-covering of $Q_2^n$ consists of $2^k$ or more
$(n-k)$-faces. If a $k$-covering $C$ of $Q_2^n$ consists of exactly
$2^{k}$ $(n-k)$-faces then $C$ is  a $k$-splitting of $Q_2^n$ into
$(n-k)$-faces. If $n-k=1$ then such splitting is equivalent to a
perfect matching in the Boolean hypercube.  A $k$-covering of
$Q_2^n$ is called an antipodal $k$-splitting if it consists of
exactly $2^{k}$ $(n-k)$-faces and it does not contain pairs of
parallel non-antipodal faces.

The concept of DP-coloring  was developed  by Dvo${\rm
\breve{r}\acute{a}}$k and Postle \cite{DP} in order to generalize
the notation of a proper coloring. In \cite{BK18} Bernshteyn and
Kostochka considered a problem to estimate the minimum number of
edges in non-$2$-DP-colorable $k$-uniform hypergraphs. The existence
of a non-$2$-DP-colorable $k$-uniform hypergraph with $e$ edges and
$n$ vertices is equivalent to the existence of a covering of $Q_2^n$
by $e$ pairs of antipodal $(n-k)$-faces. If the hypergraph has no
multiple edges then the definition of DP-coloring implies that this
covering does not contain pairs of parallel non-antipodal faces. If
$e=2^{k-1}$ then a non-$2$-DP-colorable $k$-uniform hypergraph with
$e$ edges generates an antipodal $k$-splitting and vice versa. The
connection between
  $2$-colorings of hypergraphs and coverings of the hypercube will be stated with more
  details in Section 3.

It is known (see \cite{BK18}) that for any even $k$ each $k$-uniform
hypergraph with $2^{k-1}$ edges has 2-DP-coloring. Bernshteyn and
Kostochka conjectured  that  for any odd $k\geq 3$ there exists a
non-2-DP-coloring $k$-uniform hypergraph with $2^{k-1}$ edges. The
main result of this paper is a construction of antipodal
$k$-splittings \footnote{\,Earlier, I mistakenly claimed that such
an example does not exist, i.e., each $5$-uniform hypergraph with
$16$ edges is 2-DP-colorable (see \cite{BK18}).} for any odd $k\geq
3$ and, consequently, a proof of existence
 of non-2-DP-colorable $k$-uniform
hypergraphs with $2^{k-1}$ edges.

\section{Splittings of a hypercube}

We denote a $(n-k)$-face of $Q_2^n$ by a $n$-tuple
$(a_1,a_2,\dots,a_n)$ of symbols $0,1,*$ where the symbol $*$ is
used $n-k$ times. In more details,
$(a_1,\dots,a_n)=\{(x_1,\dots,x_n) : x_i=a_i\ \mbox{if}\ a_i=0\
\mbox{or}\ a_i=1 \}$.

If $A=\{(a_1,\dots,a_n)\}$ is an antipodal $k$-splitting then
$A_\tau=\{(a_{\tau1},\dots,a_{\tau n})\}$ is an antipodal
$k$-splitting for any permutation $\tau$.  Let us agree to $*\oplus
0=*\oplus 1=*$. We define Boolean addition of $n$-tuples to act
coordinate-wise. Then for any $(n-k)$-face $a$ and any $b\in Q^n_2$
a sum $a\oplus b$ is a $(n-k)$-face of $Q^n_2$. It is clear that if
$A$ is an antipodal $k$-splitting then $A\oplus b= \{a\oplus b :
a\in A\}$ is an antipodal $k$-splitting for each $b\in Q^n_2$. We
will refer to the aboved operations as  isometries of a Boolean
hypercube.
 $A$ and
$A'$ are called equivalent antipodal $k$-splittings if $A'$ is
obtained from $A$ by an isometry.

\begin{proposition}
If there exists an antipodal $k$-covering of $Q^{n}_2$ then there
exists an antipodal $k$-covering of $Q^{n+1}_2$ with the same
cardinality.
\end{proposition}
Proof. If $A$ is an antipodal $k$-covering of $Q^{n}_2$ then $B=\{
(a_1,\dots,a_n,*) : (a_1,\dots,a_n)\in A\}$ is an antipodal
$k$-covering of $Q^{n+1}_2$. $\square$

Let $T$ be a subset of  $\{1,\dots,n\}$. An antipodal $k$-splitting
$A$ in $Q^{n}_2$ is called $t$-balanced on $T$ if any $a\in A$ has
$t$ elements $0$ and $1$ ($|T|-t$ asterisks)  in coordinates from
$T$.

\begin{proposition}\label{proDP99}
If there exist an antipodal $k_1$-splitting of $Q^{n_1}_2$ and an
antipodal $k_2$-splitting of $Q^{n_2}_2$ which is $t$-balanced on
$T$, $|T|=n_3\leq n_2$, then there exists an antipodal $(k_2
+(k_1-1)t)$-splitting of $Q^{n_2+(n_1-1)n_3}_2$.
\end{proposition}
Proof. Let $A$ be an antipodal $k_2$-splitting of $Q^{n_2}_2$ and
$B=B_0\cup B_1$ be an antipodal $k_1$-splitting of $Q^{n_1}_2$ where
sets $B_0$ and $B_1$ do not contain parallel  $(n_1-k_1)$-faces.
Consider a $(n_2-k_2)$-face $(a_1,\dots,a_{n_2})\in A$. For $i\in
T$, if $a_i=0$ we replace $a_i$ by arbitrary $b\in B_0$; if $a_i=1$
then we replace $a_i$ by arbitrary $b\in B_1$; if $a_i=*$ then we
replace $a_i$ by $\underbrace{(*,\dots,*)}_{n_1}$. So, we obtain a
set $C$ of $|A|(|B|/2)^{t}= 2^{k_2}\cdot2^{(k_1-1)t}$ tuples
 corresponding to   $m$-faces in $Q^{n_2+(n_1-1)n_3}_2$, where
$$m=n_2-k_2-(n_3-t)+ (n_3-t)n_1+ t(n_1-k_1)=n_2+(n_1-1)n_3-(k_2
+(k_1-1)t).$$ It is not difficult to verify that 1) all tuples of
$C$ are disjoint; 2) $C$ is a covering of $Q^{n_2+(n_1-1)n_3}_2$ by
counting  cardinality of $\cup C$ and, consequently, $C$ is $(k_2
+(k_1-1)t)$-splitting; 3) $C$ contains pairs of antipodal faces
because $A$ and $B$ contain pairs of antipodal faces; 4) $C$ does
not contain parallel non-antipodal faces because $A$ and $B$ do not
contain such faces. $\square$

\begin{corollary}\label{proDP999}
If there exist an antipodal $k_1$-splitting of $Q^{n_1}_2$ and an
antipodal $k_2$-splitting of $Q^{n_2}_2$ then there exists an
antipodal $k_1k_2$-splitting of $Q^{n_1n_2}_2$.
\end{corollary}

The following antipodal $3$-splitting  of $Q^{4}_2$ correspond to
 the well-known antipodal perfect matching in $Q^{4}_2$.    We
will denote it by $E_3$.

\medskip

$\begin{matrix}
*&0&0&0,\quad & *&1&1&1, \\
0&*&0&1,\quad & 1&*&1&0, \\
0&1&*&0,\quad & 1&0&*&1, \\
0&0&1&*,\quad & 1&1&0&*. \\
\end{matrix}$

\medskip

There exist two  antipodal $5$-splittings of $Q^{8}_2$.

\medskip

$\begin{matrix}
00&1*&00&**, \quad & 11&0*&11&**, \\
00&*1&10&**, \quad & 11&*0&01&**, \\
0*&01&00&**, \quad & 1*&10&11&**, \\
*0&10&01&**, \quad & *1&01&10&**, \\
01&1*&0*&0*, \quad & 10&0*&1*&1*, \\
00&*0&1*&1*, \quad & 11&*1&0*&0*, \\
0*&10&1*&0*, \quad & 1*&01&0*&1*, \\
*0&00&0*&1*, \quad & *1&11&1*&0*, \\
01&0*&*1&0*, \quad & 10&1*&*0&1*, \\
10&*1&*0&0*, \quad & 01&*0&*1&1*, \\
1*&10&*0&0*, \quad & 0*&01&*1&1*, \\
*1&00&*0&*0, \quad & *0&11&*1&*1, \\
*0&00&**&00, \quad & *1&11&**&11, \\
*0&0*&*1&01, \quad & *1&1*&*0&10, \\
*0&*1&*1&00, \quad & *1&*0&*0&11, \\
**&00&*0&01, \quad & **&11&*1&10; \\
\end{matrix}$

\bigskip

$\begin{matrix}
0*&*0&  *0&10, \quad &1*&*1 &*1&01 , \\
0*&*1&  1*&10, \quad &1*&*0 &0*&01 , \\
0*&*1&  00&*0, \quad &1*&*0 &11&*1 , \\
0*&*0&  00&0*, \quad &1*&*1 &11&1* , \\
*1&0*&  *1&00, \quad &*0&1* &*0&11 , \\
*1&1*&  1*&00, \quad &*0&0* &0*&11 , \\
*0&0*&  01&*0, \quad &*1&1* &10&*1 , \\
*1&0*&  10&0*, \quad &*0&1* &01&1* , \\
*1&*0&  *1&10, \quad &*0&*1 &*0&01 , \\
*0&*1&  1*&00, \quad &*1&*0 &0*&11 , \\
*0&*0&  11&*0, \quad &*1&*1 &00&*1 , \\
*0&*0&  10&0*, \quad &*1&*1 &01&1* , \\
1*&1*&  *0&10, \quad &0*&0* &*1&01 , \\
1*&1*&  0*&00, \quad &0*&0* &1*&11 , \\
1*&0*&  00&*0, \quad &0*&1* &11&*1 , \\
0*&1*&  01&0*, \quad &1*&0* &10&1* . \\
 \end{matrix}$

\medskip

Note that any isometry of the hypercube exchange only the order of
columns and symbols $0$ and $1$ in any fixed column. Consequently,
the above two antipodal $5$-splittings are nonequivalent.  The last
antipodal $5$-splitting is $1$-balanced on the sets $\{1,2\}$ and
$\{3,4\}$. We will denote it by $E_5$.

\begin{theorem}\label{thDPcor1}
There exists an antipodal $k$-splitting  for every odd $k\geq 3$.
\end{theorem}
Proof.
 Let us use antipodal $5$-splitting $E_5$ and $3$-splitting $E_3$ in the construction from
 Proposition \ref{proDP99}
 with $k_1=3, k_2=5$,
 $n_1=4, n_2=8$, $T = \{1,2\}$, $t=1$. We obtain an antipodal  $7$-splittings
 $E_7$. Replacing  $E_3$ by $E_5$ in this construction, we obtain an antipodal  $9$-splittings
 $E_9$. It is clear that we can construct an antipodal $(2s+1)$-splittings
 $E_{2s+1}$ from $(2s-3)$-splittings
 $E_{2s-3}$ and $1$-balanced $5$-splittings $E_5$ by Proposition
 \ref{proDP99}. Consequently, we prove the theorem by induction.
 $\square$

The proof of the following  statement can be found in  \cite{BK18}.
Here it is rewritten  in  notations of this article.

\begin{proposition}[\cite{BK18}]
If $k$ is even then an antipodal $k$-splitting of $Q^n_2$ does not
exist.
\end{proposition}
Proof. Let $A$ be an antipodal $k$-splitting and $k$ is even. Let us
consider an $(n-k)$-face $a\in A$, the $(n-k)$-face $\overline{a}\in
A$  antipodal to $a$ and a $k$-face $a^\bot$ orthogonal (dual) to
$a$, i.e., positions of asterisks in $a$ and $a^\bot$ are
complementary. By the definitions, we obtain that $x=a\cap a^\bot$
and $\widetilde{x}=\overline{a}\cap a^\bot$ are antipodal vertices
within the face $a^\bot$. For example, $a=(0,1,1,0,*,*)$,
$a^\bot=(*,*,*,*,1,0)$, $x=(0,1,1,0,1,0)$,
$\widetilde{x}=(1,0,0,1,1,0)$. The vertices $x$ and $\widetilde{x}$
have the same parity because $k$ is even. But for all other $b\in A$
we obtain that $b\cap a^\bot$ has the same number of even-weighted
and odd-weighted vertices. Since $A$ is a splitting, the set
$\{b\cap a^\bot : b\in A\}$ is a splitting of $a^\bot$ as well.
Because the numbers of even-weighted and odd-weighted vertices in
$a^\bot$ are equal, we have a contradiction. $\square$

\begin{proposition}\label{DPdop222}
For any $k$-splitting $A$ of $Q^n_2$  ($n>k>0$) and for any
direction of faces the number of $(n-k)$-faces of this direction in
$A$ is even.
\end{proposition}
Proof.
 Suppose $a\in A$ and $A$ contains $m$
$(n-k)$-faces of the same direction as $a$. Consider a face
$a^\bot$. If $b\in A$ has the same direction as $a$ then $|b\cap
a^\bot|=1$, otherwise the number $|b\cap a^\bot|$ is even. Since
$|a^\bot|=2^k=\sum_{b\in A} |b\cap a^\bot|$ and all terms except $m$
are even, $m$ is even. $\square$

A splitting of hypercube into $(n-k)$-faces is a special case of
A-designs. In \cite{P14} there are given constructions of A-design
with additional properties, for example, with no adjacent parallel
faces. An antipodal $k$-splitting does not exist if $k$ even.
Nevertheless we can find a $k$-splitting of $Q^n_2$ with a pair of
parallel $(n-k)$-faces on almost maximal  distance.

\begin{proposition}\label{DPdop44}
For  even $k\geq 4$ there exits a $k$-splitting of $Q^n_2$ with the
distance $k-1$ between every pair of $(n-k)$-faces of the same
direction.
\end{proposition}
Proof. As mentioned above, for any permutation $\tau$ and any
$k$-splitting $A$ a set $A_\tau$ is a $k$-splitting. Consider
$k$-splitting $(E_5)_\tau$ for $\tau=
\begin{pmatrix}
1234 & 5678 \\
5678 & 1234\\
\end{pmatrix}$.
Since the first four coordinates of faces from $E_5$ contain two
symbols $*$ and the last four coordinates of faces from $E_5$
contain one symbol $*$, we conclude that  $(E_5)_\tau$ does not
contain faces of the same direction as $E_5$. Denote a set
$E_50=\{a0 : a\in E_5\}$. So a distance between parallel faces in a
$6$-splitting $E_50\cup (E_5)_\tau1$ is maximal minus $1$, i.e., it
is equal to $5$. By the same way, we can prove that $E_{2s+1}**$ and
$(E_{2s+1}**)_\sigma$, where $\sigma$ permute two last pairs of
coordinates, do not contain faces of the same direction. Then
$E_{2s+1}**0\cup(E_{2s+1}**)_\sigma1$ is the required
$(2s+2)$-splitting. $\square$

Now we will find a minimal  dimension of a hypercube  containing a
$k$-splitting  with at most two $(n-k)$-faces of any fixed
direction.

\begin{proposition}\label{DPdop22}
There exits a $k$-splitting of $Q^n_2$ with at most two
$(n-k)$-faces of each direction if $n-2k+2\geq 0$, $0<k<n$.
\end{proposition}
Proof. By induction on $n$.  Any $1$-splitting of $Q^n_2$ consists
of two parallel faces. For $n=3, k=2$ it is easy to verify this
statement directly. The case $n=4$, $k=3$ corresponds to $E_3$. Let
$B$ be a $k$-splitting of $Q^{n}_2$ with at most two $(n-k)$-faces
of each direction. By Proposition \ref{DPdop222} $B$ contains two or
zero faces of each direction. Let $B=B_0\cup B_1$ where sets $B_0$
and $B_1$ do not contain parallel $(n-k)$-faces. Consider the set
$A=\{b0*,b1* : b\in B_0\}\cup \{b*1,b*0 : b\in B_1\}$. By the
construction, the set $A$ is a $(k+1)$-splitting of $Q^{n+2}_2$ with
at most two $(n-k+1)$-faces of each direction. Besides, the set
$C=B*$ is a $k$-splitting of $Q^{n+1}_2$. $\square$

\section{2-DP-colorings}

Let $G$ be an $r$-uniform hypergraph on $n$ vertices. For each $e\in
E(G)$ we consider two antipodal 2-colorings $\varphi_e:e\rightarrow
\{0,1\}$ and $\overline{\varphi_e}=\varphi_e\oplus 1$. Let $\Phi$ be
a collection of $\varphi_e$, $e\in E(G)$. We say that a 2-coloring
$f: V(G)\rightarrow \{0,1\}$ avoids $\Phi$ if $f|_e\neq \varphi_e$
and $f|_e\neq \overline{\varphi}_e$ for each $e\in E(G)$.

A hypergraph $G$ is called a proper 2-colorable if there exists a
$2$-coloring $f$ avoiding $\Phi_0$, where $\Phi_0$ consists of
constant maps. A hypergraph $G$ is called {2-DP-colorable} if for
every  $\Phi$ there exists a 2-coloring $f$ avoiding $\Phi$.

A $2$-coloring $f$ of a $k$-uniform hypergraph on $n$ vertices is
one-to-one corresponding to an $n$-tuple over alphabet $\{0,1\}$
($f\in Q^n_2$). Each $k$-hyperedge  correspond to  $(n-k)$-faces of
$Q^n_2$ of some direction. For example,  a $k$-hyperedge consisting
of $i_1$th,...,$i_k$th vertices corresponds to faces
$(*,\dots,*,a_{i_1},*,\dots,*,a_{i_2},*,\dots,a_{i_k},*\dots,*)$
where $a_{i_j}\in \{0,1\}$. The set $\{a_{i_j}\} $ corresponds to
some coloring of vertices from the hyperedge. A $2$-coloring $f$
avoids $\varphi_e = (*,\dots,1,\dots, *,\dots, 0,\dots,*)$ iff
$f\not \in \varphi_e$. A 2-coloring $f$ avoids $\Phi$ if $f\not \in
\varphi_e\cup \overline{\varphi}_e$ for each $\varphi_e\in \Phi$.

Consider a table of size $n\times \ell$,   where every column
corresponds to a $(n-k)$-face of an antipodal covering of $Q^n_2$.
Let us replace in the table symbols $*$ by $0$ and other symbols by
$1$. By the definition, the resulting table is the incidence matrix
of a non-$2$-DP-colorable $k$-uniform hypergraph with $\ell$ edges.
Consequently, we have the following statement.

\begin{proposition}\label{proDP103}
A  $k$-uniform hypergraph with $\ell$ edges and $n$ vertices is
non-$2$-DP-colorable if and only if its incidence matrix corresponds
to
  a $k$-covering of $Q_2^n$ by $\ell$ pairs
of antipodal $(n-k)$-faces.
\end{proposition}

Moreover, Proposition \ref{proDP103} implies the following
statement.

\begin{corollary}\label{proDP102}
There exists  a non-$2$-DP-colorable $k$-uniform hypergraph with
$2^{k-1}$ edges  if and only if there exists  an antipodal
$k$-splitting of $Q^n_2$.
\end{corollary}

 A non-$2$-DP-colorable $3$-uniform hypergraph with
$4$ edges that corresponding to the antipodal $3$-splitting $E_3$ is
presented  in \cite{BK18}. By Theorem \ref{thDPcor1} and Corollary
\ref{proDP102},  we obtain the following statement.

\begin{corollary}\label{DPcor0}
If $k\geq 3$ is odd there exists a  non-2-DP-colorable $k$-uniform
hypergraph  with $2^{k-1}$ edges.
\end{corollary}

Since a union of  at most $2\ell$ $(n-k)$-faces contains $\ell
2^{n-k+1}$ vertices, we obtain the following corollary.

\begin{corollary}[\cite{BK18}]
Every $k$-uniform hypergraph with  $\ell <2^{k-1}$ edges is
$2$-DP-colorable.
\end{corollary}

By definition a proper coloring of hypergraph is not contained
monochromatic edges. Consequently,
 a non-$2$-proper colorable hypergraph corresponds to a covering of
the hypercube consisting of faces which contain vertices
$\overline{0}$ or $\overline{1}$. Therefore, each $k$-uniform
hypergraph with $2^{k-1}$ edges is proper $2$-colorable. Moreover,
by similar arguments we obtain that every $k$-uniform hypergraph
with $s^{k-1}$ or less edges is proper $s$-colorable.  A better
bound for the case of proper colorings is known. Cherkashin and
Kozik \cite{ChK} Radhakrishnan and Srinivasan \cite{RS} (for $s=2$)
showed that any $k$-uniform hypergraph with $c(s)(\frac{k}{\ln
k})^\frac{s-1}{s}s^{k-1}$ or less edges is proper $s$-colorable,
where $c(s)>0$ does not depend on $k$ ($k$ is large enough). A
survey of results on  proper colorings of hypergraphs and related
problems were found in \cite{RSh}. A some Brooks' type theorem for
DP-colorings of hypergraphs  is proved in \cite{Schweser}.

\section{Acknowledgments}

The author is grateful to S.Avgustinovich, A.Kostochka and
A.Taranenko for their attention to this work and  useful
discussions.


\begin{thebibliography}{9}








\bibitem {BK18}
Bernshteyn A., Kostochka A. ``DP-colorings of hypergraphs, European
Journal of Combinatorics'', V. 78  (2019), P.~134--146.


\bibitem {ChK}
  Cherkashin D.,  Kozik J. ``A note on random greedy coloring of uniform hypergraphs'', Random
Structures and Algorithms, V. 47:3 (2015), P.~407-416.


\bibitem {DP}
Dvo${\rm \breve{r}\acute{a}}$k Z.,   Postle L. ``Correspondence
coloring and its application to list-coloring planar graphs without
cycles of lengths 4 to 8'', Journal of Combinatorial Theory. Series
B, V. 129 (2018),
P.~38--54. 











\bibitem {RS}
Radhakrishnan J., Srinivasan  A. ``Improved bounds and algorithms
for hypergraph two-coloring'', Random Structures and Algorithms, V.
16 (2000), P.~4-32.


\bibitem {Schweser}
Schweser T. ``DP-degree colorable hypergraphs'', Theoretical
Computer Science, V. 796 (2019), P.~196-206.





\bibitem {P14}
Potapov V.N. ``On the multidimensional permanent and $q$-ary
designs'',
 Siberian Electronic Mathematical Reports.,  V. 11 (2014),
P.~451--456.


\bibitem {RSh}
 Raigorodskii A.M.,  Shabanov D.A., ``The Erd\"os-Hajnal problem of
hypergraph colouring, its generalizations, and related problems'',
 Russian
Math. Surveys, V. 66:5 (2011), P.~933--1002.




\end{thebibliography}
\end{document}